\documentclass[runningheads]{PMS2026}

\usepackage[ansinew]{inputenc}
\usepackage[T1]{fontenc}
\usepackage[top=3cm,bottom=3.5cm,left=4cm,right=3.8cm,dvips]{geometry}
\usepackage{nopageno}

\usepackage{index}
\usepackage{fancyhdr}

\usepackage{indentfirst}
\usepackage{sectsty}
\allsectionsfont{\usefont{OT1}{phv}{bc}{n}\selectfont \normalsize \bf }
\usepackage{graphicx}
\usepackage{subfigure}
\usepackage{url}
\usepackage{harvard}
\usepackage{amssymb}
\usepackage{epsfig}

\usepackage{amsmath, tikz}
\usetikzlibrary{calc, positioning}

\makeindex

\begin{document}
\pagenumbering{gobble}
%\pagenumbering{roman}
\pagestyle{headings}
\mainmatter

%\tableofcontents

\title{Solving a large oral examination timetabling problem using a multidimensional knapsack MILP formulation}

\author{ Cyrille Briand\inst{1} and Jean-Pierre Belaud\inst{2} }

\institute{LAAS-CNRS, University of Toulouse, France\\
\email{ cyrille.briand@laas.fr}
\and
LGC, University of Toulouse, CNRS, Toulouse, France\\
\email{jeanpierre.belaud@toulouse-inp.fr}
}

\maketitle
% Please, write again your names here under, as you want it to appear in the author's index.
\index{Cyrille Briand}\index{Jean-Pierre Belaud}

\noindent
\textbf{Keywords:} Examination timetabling, multidimensional knapsack , mixed integer linear programming.

\section {Introduction}
France's prestigious engineering schools, known as \emph{Grandes \'{E}coles},  train students from the third year of bachelor's degree (L3) to the second year of a master's degree (M2). Admission is determined by performance in  major groups of written examinations and, if successful, oral examinations. Each group targets a distinct student profile and range of schools. The largest academic track, MP (Mathematics and Physics), historically accounts for the majority of candidates, with approximately 6,500 to 7,000 students taking the exams annually.The largest academic track, MP (Mathematics and Physics), historically accounts for the majority of candidates, with approximately 6,500 to 7,000 students sitting the exams each year. The newer track, MPI (Mathematics, Physics, and Computer Science) has a more modest contingent, around 1,000 candidates annually. While both tracks share the same written exam bank, their specific oral exams reflect the stronger emphasis on computer science for the MPI track.

The organization of oral examinations is entrusted to the Service des Concours Communs Polytechniques (SCCP), hosted by Toulouse INP University.. The SCCP is responsible for operating the oral examinations of CentraleSup\'{e}lec, Mines-Ponts, and Concours Commun INP (CC INP), which provide access to a wide range of schools (approximately 80 schools with 5,700 available places). The organization is complex, as a single student may be eligible for several oral examination groups, each with its own constraints, requiring careful coordination of the individual exam timetables.
 
The oral examinations consist of offering each successfully admitted candidate a schedule of oral exams, i.e., a specific time slot for each exam for which they are admissible. The timeline is very tight, as schedules must be built within only a few days following the release of the admissibility results. During each scheduled slot, a candidate may take several tests (e.g., practical work, modern language oral exams). Therefore, the duration of a slot depends on the type of the exam. This is a large scale schedule problem as there are more than 7,000 candidates to accommodate over a 4-weeks time horizon. It is also complex as each candidate can be eligible for several types of oral exam, which requires to avoid overlapping slots. Furthermore, the specific modern language chosen by the candidate must be carefully considered, as not all languages are available in every slot. In addition, for a given exam type, there are several possible slots, and each is characterized by its own  constraints:   time length,  type of  slot (MP, MPI, or both),  maximum number of candidates, modern languages available,  maximum number of candidates for each language. Finally, since the oral examinations take place in Paris, the proposed schedule must account for the candidate's place of origin. For instance, since accommodation in Paris is expensive, the slots offered to candidates from outside Paris must avoid long idle times. For non-French candidates, the time required to obtain a visa must be accounted for, meaning that only late slots, scheduled towards the end of the four-week period, should be proposed.

This paper tackles this peculiar Oral Examination Timetabling Problem (OETP). It first shows how this problem can be advantageously modeled as a Multidimensional Knapsack Problem (MKP), and the pros and cons of such an approach are discussed. The paper also provides an initial set of experiments on the real 2025 SCCP MP/MPI problem instance, utilizing an open-source MILP solver. These initial results demonstrate the effectiveness of the approach. The remainder of the paper is structured as follows. Section~\ref{sec:lit_review} presents a brief literature review of the OETP. Section~\ref{sec:mod} details our BMP-oriented modeling, which is formulated as a mixed integer linear program. Section~\ref{sec:exp} details the experimentation made. Conclusion and perspectives are drawn in Section~\ref{sec:ccl}.

\section{Brief literature Review}
\label{sec:lit_review}
\sloppy

The Oral Examination Timetabling Problem (OETP) can be viewed as a special case of the University Examination Timetabling Problem (UETP), which has been extensively studied over the past decades due to its inherent computational complexity. UETP is defined as the assignment of a set of exams into a limited number of time slots and rooms while satisfying a variety of hard and soft constraints, making it a classic NP-hard combinatorial optimization problem. Seminal work by \cite{carter96} established the foundational methodology for UETP, primarily modeling the problem as a Graph Coloring Problem. In this paradigm, exams are represented as vertices and conflicts (shared students) as edges; the objective is to assign colors (time slots) such that no adjacent vertices share the same color. Following these early graph-theoretic formulations, the literature has expanded to include a wide variety of solution techniques. As detailed in comprehensive surveys such as the one by Qu et al. \cite{qu09}, approaches have evolved from constructive heuristics to sophisticated metaheuristics, including evolutionary algorithms, simulated annealing, and tabu search. Moreover, as discussed in \cite{Muller16} based on a real large-scale UETP, Constraint Programming (CP) is a powerful approach capable of conveniently handling most complex constraints, yet scalability remains a challenge for large instances.

Most existing studies address the room assignment problem using matching algorithms or network flow models, typically after the time slots have been fixed. Few studies, however, explore the simultaneous optimization of time and resources using a packing perspective, although assigning exams to rooms with heterogeneous capacities shares structural properties with the Knapsack Problem. The Multidimensional Knapsack Problem (MKP) is a well-known variation of the Knapsack Problem where items (exams) consume multiple resources (e.g., seats, invigilators, equipment) simultaneously. Given the increasing complexity of modern university requirements, where constraints are not merely about "fitting" students into rooms but about optimizing multi-attribute resources, an MKP-oriented approach could be relevant. Therefore, this paper investigates such an approach, tackling the OETP using an MKP-based Mixed Integer Linear Programming (MILP) formulation.

\section{Modelling}
\label{sec:mod}

To model the OETP, one could employ a direct method involving the creation of an assignment variable for each candidate-slot combination. While this method offers the substantial advantage of ensuring comprehensive consideration of all potential assignments, it suffers from significant drawbacks as it leads to a combinatorial explosion due to the creation of an excessive number of variables and constraints, drastically increasing the overall model complexity. Specifically, integrating intricate operational constraints, such as non-overlapping slots or candidate origin constraints, is not trivial. 

Alternatively, our approach considers only a predetermined set of existing timetables, potentially derived from historical data, and assigns one of these pre-validated schedules to each candidate. The primary benefit of this strategy is its ability to severely restrict the total number of variables and constraints, thereby simplifying the model. It also naturally supports a human-in-the-loop validation process, as schedules can be specifically designed to accommodate candidates. However, its main limitation is the requirement for a separate mechanism to dynamically construct new consistent schedules.

Formally, the OETP can be defined by the tuple $(\mathcal{C},\mathcal{P},\mathcal{R},\{\mathcal{P}_i\},\{\mathcal{P}_k\})$. $\mathcal{C}$ is the set of \emph{candidates}, $\mathcal{P}$ denotes the set of \emph{schedules} (or timetables), and $\mathcal{R}$ is the set of resources, each resource $k\in \mathcal{R}$ having a capacity $b_k$. $\mathcal{P}_i \subseteq \mathcal{P}$ is the subset of compatible schedules for candidate $i$. Symmetrically, $\mathcal{P}_k \subseteq \mathcal{P}$ is the subset of schedules that utilize resource $k$.

\begin{figure}[htbp]
\begin{center}
\includegraphics[scale=0.3]{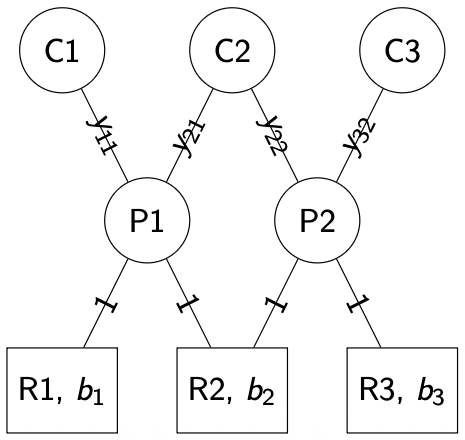}
\caption{Illustration of the formulation}
\label{fig}
\end{center}
\end{figure}

The objective is to maximize the total number of assigned candidates subject to capacity constraints. This problem is mathematically formulated as a Multidimensional Knapsack Problem (MKP) below. The notation is defined as follows (see Figure~\ref{fig} for illustration). $y_{ij} \in \{0,1\}$ are binary decision variables. $y_{ij}=1$ if candidate $i$ is assigned to schedule $j \in \mathcal{P}_i$, and $0$ otherwise. The objective function \eqref{obj-func} expresses the goal of maximizing the total number of assigned candidates. Constraints \eqref{res-cons} are the resource capacity constraints, ensuring that the total usage of resource $k$ by the selected schedules does not exceed the capacity $b_k$. Finally, constraints \eqref{cand-cons} are the assignment constraints, ensuring that each candidate $i$ is assigned to at most one schedule $j \in \mathcal{P}_i$.

\begin{align}
\max \quad & \sum_{i\in \mathcal{C}} \sum_{j\in \mathcal{P}} y_{ij} \label{obj-func}\\
\text{subject to} \quad & \sum_{i\in \mathcal{C}}\sum_{j \in \mathcal{P}_k} y_{ij} \leq b_k, \quad \forall k \in \mathcal{R} \label{res-cons}\\
\quad & \sum_{j \in \mathcal{P}_i} y_{ij} \leq 1, \quad \forall i \in \mathcal{C} \label{cand-cons}\\
& y_{ij} \in \{0,1\}, \quad \forall (i,j) \in \mathcal{C} \times \mathcal{P}
\end{align}

\section{Experiments}
\label{sec:exp}
This experiment evaluates the performance of SCIP (Solving Constraint Integer Programs) on a large-scale assignment problem derived from the 2025 SCCP examination process. SCIP is a state-of-the-art framework for integer and linear optimization. It is freely available for academic use and distributed with its full source code. Its Apache 2.0 license allows unrestricted use, including in commercial contexts. 

The studied MP/MPI instance involves 7,804 candidates, each characterized by a highly heterogeneous profile. The overall search space comprises 7,759 feasible schedules (or timetables). A dedicated greedy heuristic was implemented to explore all admissible assignments, derived from the 1,123,321 candidate-schedule compatibility relations. Each assignment simultaneously utilizes up to four resources from a pool of 103 available resources, each with its specific capacity. The greedy procedure produces an initial feasible solution, assigning 7,759 candidates (out of the 7,804 total). This solution utilized 606 distinct schedules. This solution was injected into SCIP as a warm start, along with the 606 schedules actually utilized by the greedy solution. A set of 100 still unsaturated plannings were also provided to SCIP, enabling the solver to focus its search on the most promising region of the solution space, under a time limit of 20 minutes. SCIP successfully computed an \emph{optimal} solution accommodating 7,796 candidates among 7,804, thus demonstrating the benefit of combining a fast constructive heuristic with a powerful MILP solver. Assigning more candidates would require increasing the pool of available schedules (or to relax unsatisfiable constraints).

%This experiment represents an initial proof-of-concept. Future work will investigate more sophisticated modeling and relaxation strategies, especially considering that, in real operational settings, some constraints may be deliberately softened (or relaxed) to ensure that all candidates can be assigned.

\section{Conclusion and Future Work}
\label{sec:ccl}
Formulating the OETP as an MKP provides a rigorous mathematical structure, enabling the effective management of the complexity associated with capacity and resource constraints. This approach is particularly relevant in planning contexts where multiple resource dimensions compete for a set of discrete items (schedules). By framing the selection of schedules under these multiple resource constraints, the MKP model not only facilitates the exploration of the solution space but also offers theoretical guarantees for the application of exact optimization algorithms and specialized heuristics, thus validating its relevance for solving complex application cases.

To overcome the computational limitations inherent in large-scale instances, a promising direction is the adoption of a Branch-and-Price approach. This decomposition allows the Restricted Master Problem to be treated as an MKP instance, while the resulting pricing subproblem focuses on efficiently generating new columns (improving schedules) with a positive reduced cost.  Once a set of interesting columns is determined, future directions could expand to integrating auxiliary objective functions, such as establishing a balanced parity across the involved academic tracks, thereby better aligning the final solution with operational needs.

%\printindex
\end{document}